\documentclass[submission,copyright,creativecommons]{eptcs}
\providecommand{\event}{ThEdu20} 
\usepackage{breakurl}             
\usepackage{underscore}           
\usepackage[utf8x]{inputenc}
\usepackage{amsmath}
\usepackage{tcolorbox}
\usepackage{amsfonts}
\usepackage{url}
\usepackage{tikz}

\title{Number Theory and Axiomatic Geometry in the Diproche System}
\author{Merlin Carl
\institute{Institut f\"ur mathematische, naturwissenschaftliche und technische Bildung, Abteilung f\"ur Mathematik und ihre Didaktik\\
Europa-Universit\"at Flensburg\\
Flensburg, Germany}
\email{merlin.carl@uni-flensburg.de}}

\def\titlerunning{Number Theory and Geometry in Diproche}
\def\authorrunning{M. Carl}
\begin{document}
\maketitle

\begin{abstract}
Diproche (``Didactical Proof Checking'') is an automatic system for supporting the acquistion of elementary proving skills in the initial 
 phase of university education in mathematics. A key feature of Diproche~-- which is designed by the example of the Naproche system developed by M. Cramer and others (see, e.g.,  \cite{CKKS})~-- is an automated proof checker for proofs written in a controlled fragment of natural language specifically designed to capture the language of beginners' proving exercises in mathematics. 
 Both the accepted language and proof methods depend on the didactical and mathematical context and vary with the level of education and the topic proposed. An overall presentation of the system in general was given in \cite{CK}. Here, we briefly recall the basic architecture of Diproche and then focus on explaining key features and the working principles of Diproche in the sample topics of elementary number theory and axiomatic geometry. 
\end{abstract}

\section{Introduction}

It is well-known to anyone teaching introductory classes in mathematics at the university level that understanding the concept of mathematical proof and learning how to prove is a considerable difficulty for the beginner. This ability is learned through practice; for this reason, regular exercises form an integral part of mathematical lectures. In order for this to become effective, however, feedback on the student's performance is required. If this is provided by human correctors, as it is usually the case, the time span between writing a solution and receiving the feedback is quite long (typically about a week); moreover, the feedback only comes after an exercise is finished and thus does not help while working on an exercise, e.g., by modifying a failed approach to a more successful one, or attempting another approach, or filling in details, dealing with an overlooked case, or even just improving the presentation of an argument.

The goal of the Diproche system, first introduced in Carl and Krapf \cite{CK}\footnote{Note, however, that several system components have considerably developed since \cite{CK} and that the description given there does not in every detail apply to the new version. For example, the goal tracing now deals with multiple goal candidates for each proof line rather than a single goal. Moreover, the language check was not yet implemented when \cite{CK} was submitted.}, is to provide a tool for teaching how to prove, which includes properly expressing the proof in natural language. Users are presented with a proving exercise and can enter their solution in a text window in (a controlled fragment of) natural language.\footnote{As Diproche is developed for German students, the input language is currently German. It would be unproblematic to adapt the natural language processing components to work with other languages, such as English.} The text is then translated into a formal representation format. From this, a series of proving tasks is extracted, which are given to an automated theorem prover (ATP). The ATP is carefully hand-crafted in order to accept those and only those steps that can be regarded as elementary for the respective topic and degree of education. The user is then given feedback informing her or him about (i) whether all of the steps could be verified and if not, which ones failed (ii) whether the non-verifiable steps could be explained as instances of known formal fallacies (using the ``anti-ATP'', see below), (iii) whether the announced goal of the proof was reached (if a goal was declared), (iv) whether all variables were introduced and used correctly with respect to types. 
Additionally, users can request various kinds of hints; there are also two sub-programs for learning the use of logical formalism. All of these components are briefly explained below.



Here, we will be mostly concerned with the proof-checking component and describe this for the topics of elementary number theory and axiomatic geometry, both of which are implemented in the current version of Diproche. 

A general observation behind Diproche is that, in order to provide a natural environment for expressing mathematics, a foundational perspective striving to come up with a single uniform framework for expressing all of mathematics has to be given up in favour of a variety of contexts, each with its own linguistical and logical peculiarities. In Diproche, this is realized by the so-called ``playing fields''. A playing field consists of a specific grammar and vocabulary, along with basic notational conventions and inference rules adapted to a certain mathematical topic. Proving exercises then always take place in a certain ``playing field''. Among other advantages, this allows us to use the same notation differently in different contexts and to impose variying requirements of strictness concerning logical notation for different areas. In the current version, the following ``playing fields'' have been implemented:

\begin{itemize}
\item Propositional Logic (For exercises like: ``Show that $A\rightarrow(B\rightarrow A)$ for all propositions $A$, $B$'')
\item Boolean Set Theory (For exercises like: ``Show that $A\cap(B\cup C)=(A\cap B)\cup(A\cap C)$ for all sets $A$, $B$, $C$'')
\item Functions and Relations (For exercises like: ``Show that, if $f:B\rightarrow C$ and $g:A\rightarrow B$ are injective, then so is $f\circ g$'')
\item Elementary Number Theory (For exercises like: ``Show that, for every integer $n$, $n^2-n$ is even'')
\item Induction (For exercises like: ``Show that, for all natural numbers $n$, we have $2^{n+5}>(n+5)^2$'')
\item Axiomatic Geometry (see below for example exercises)
\item Group theory (For exercises like: ``Show that, if $(G,\cdot)$ is a group with neutral element $e$ and $a,b\in G$ with $a\cdot b=e$, then $b\cdot a=e$'')
\end{itemize}


For each playing field, one can now set up exercises. Formally, an exercise is an $8$-tuple 

\begin{center}
(Id, Nat, Form, Diff, Assmpts, Decls, PF, Hints), 
\end{center}

\noindent
where Id is the identifier of the problem, Nat is the natural language formulation displayed to the user, Form is a formalization of the goal of the exercise in the underlying Prolog format, Diff is the degree of difficulty, i.e., the set of allowed inference rules, Assmpts is a list of statements that may be used in the proof (since they are, e.g., known from the lecture), Decls is a list of declarations of variables, functions and predicates that may be used in writing the solution, PF is the identifier of the playing field to which the problem belongs and Hints is a list of hints for the user written in natural language that can be displayed on the users' request.

In this article, we will present our work so far, along with the results, on implementing the ``playing fields'' on elementary number theory and axiomatic geometry, which are the most advanced playing fields in the current version, both linguistically with respect to the size of the vocabulary and the complexity of the grammar and with respect to the number of ATP-rules required to allow for typical solutions to be handled in the way that they should. Since Diproche works in German, sample texts will be German, with the only adaptation that we replace the input format accepted by the Diproche interface with {\LaTeX} for the sake of this article. In some cases, we provide English translations for these text examples.\footnote{In order to give the reader a concrete impression of the actual Diproche language, we also give the original German texts as accepted by Diproche.} It should be noted that these can \emph{not} be accepted by the current (German) Diproche system, although it would not be hard to implement an English version of Diproche, for which the accepted texts would be very similar to these translations.

\section{Natural Language Proof Checking} 

In this section, we go into some detail with respect to the proof checking function of Diproche. The guiding idea of the architecture is the same as for the Naproche system (see, e.g., \cite{CKKS}), but the details are different and the system was implemented from scratch. The main reason is that, although quite impressive in power, Naproche is not well-adapted to didactical uses (which also was not its purpose). 


The following diagram gives the overall structure of the Diproche system.\footnote{For the sake of greater clarity, the type-checking component has been omitted.}

\begin{tikzpicture}[scale=0.7]

\draw (10,3)--(18,3)--(18,4)--(10,4)--(10,3); 
\node at (14,3.5) {Interface};
\draw (13,3)--(13,4); 
\node at (11.5,3.5) {Input};
\draw (15,3)--(15,4); 
\node at (16.5,3.5) {Output};

\draw (11,0)--(17,0)--(17,1)--(11,1)--(11,0); 
\node at (14,0.5) {Language Check};

\draw [->] (12,3)--(12,1); 
\draw [<-] (16,3)--(16,1); 

\draw (19,0)--(25,0)--(25,1)--(19,1)--(19,0); 
\node at (22,0.5) {Formula Parsing};
\draw [<->] (19,0)--(17,-1); 

\draw (19,-2)--(25,-2)--(25,-1)--(19,-1)--(19,-2); 
\node at (22,-1.5) {Text Parsing};
\draw [<->] (19,-2)--(17,-3); 

\draw (11,-1)--(17,-1)--(17,-2)--(11,-2)--(11,-1); 
\node at (14,-1.5) {Preprocessing};

\draw [<-] (14,-1)--(14,0); 

\draw (11,-3)--(17,-3)--(17,-4)--(11,-4)--(11,-3); 
\node at (14,-3.5) {Annotation};

\draw [->] (14,-2)--(14,-3); 

\draw (11,-5)--(17,-5)--(17,-6)--(11,-6)--(11,-5); 
\node at (14,-5.5) {Text structure};

\draw [->] (14,-4)--(14,-5); 

\draw [->] (17,-3.5)--(17.5,-3.5)--(17.5,-7.5)--(17,-7.5); 
\draw [->] (17,-3.5)--(20,-3.5)--(20,-9);  
\draw [->] (17,-5.5)--(19,-5.5)--(19,-9); 

\draw (11,-7)--(17,-7)--(17,-8)--(11,-8)--(11,-7); 
\node at (14,-7.5) {Generating ATP-Tasks};

\draw [->] (14,-6)--(14,-7); 

\draw (11,-9)--(17,-9)--(17,-10)--(11,-10)--(11,-9); 
\node at (14,-9.5) {ATP};

\draw [->] (14,-8)--(14,-9); 
\draw [->] (14,-10)--(14,-13); 

\draw (19,-9)--(25,-9)--(25,-10)--(19,-10)--(19,-9); 
\node at (22,-9.5) {Goal Check};

\draw [->] (20,-10)--(20,-13); 

\draw (13,-13)--(21,-13)--(21,-14)--(13,-14)--(13,-13); 
\node at (17,-13.5) {Feedback};

\draw [->] (21,-13.5)--(26,-13.5)--(26,3.5)--(18,3.5); 

\end{tikzpicture}

The function of these components will be explained below. We mention here that the ``ATP'' module actually consists of various submodules, one for each of the topics of propositional logic, Boolean set theory, functions and relations, elementary number theory and axiomatic geometry, with an additional module for algebraic term manipulations. 

Additionally, for each exercise, a set of assumptions and a set of available inference rules can be specified. Given the evolving nature of inferential abilities in mathematical education, this is necessary: Mathematics has the habit of turning results into methods, and sometimes, the proof goal of one exercise should be available as an inferential step in the next. For example, while the compatibility of parallelism and perpendicularity may be an exercise in an early stage of  learning axiomatic geometry, one should later be able to simply use it in proofs, even without mentioning it. As Descartes' put it (see \cite{De}, second part, p. $12$): ``each truth that I found being a rule that later helped me to find others.''.
Thus, one cannot use \emph{one} set of ATP-rules, even for a fixed topic; instead, the set of available ATP-rules has to develop 
while the student acquires new skills. In Diproche, this is made possible by specifying a ``difficulty degree'' for an exercise (see above). 


\subsection{Processing Example}

Let us have a look at a simple text example to see how it is processed by the different components introduced above.\footnote{A similar discussion with a different example text was given in Carl and Krapf \cite{CK}. Since \cite{CK} was in German, we give here a brief explanation for the sake of English readers.}

The user might type the following input into a text window:\footnote{Here, a remark is in order to clarify the relation of using a controlled natural language and the input accepted by the system: While most functions of the system (such as logical checking, goal tracing or type-checking) only work on texts written in the controlled natural language, users are free to enter whatever they want, including random strings. On this arbitrary input string, a first check for unknown symbols, unknown words, non well-formed formulas and syntactically non-processable sentences is performed. If one of these difficulties occurs, the user is provided with the according feedback and no further processing takes place.}

\begin{tcolorbox}
Es seien $g$, $h$, $l$ Geraden.
Angenommen, wir haben $g\parallel l$ und $l$ ist ausserdem orthogonal zu $h$. 
Dann ist auch $g$ orthogonal zu $h$.\footnote{Let $g$, $h$, $l$ be lines. Suppose we have $g\parallel l$ and furthermore that $l$ is perpendicular to $h$. Then $g$ is also perpendicular to $h$.}
\end{tcolorbox}

The language check runs through the text to determine whether there are any unknown symbols or words. If so, the processing would be stopped and an error message returned. In this text, this is not the case, so the system continues with preprocessing.

The preprocessing turns this string into a list of sentences, each of which is in turn represented as a list of words. Moreover, formal expressions are identified as either terms or formulas, converted into an internal list representation and assigned with their type. Thus, we get the following output from the preprocessing:

\bigskip

[[es,seien,g,h,l,geraden],

[angenommen,wir,haben,[fml,[g,parallel,l]],und,l,ist,ausserdem,orthogonal,zu,h],

[dann,ist,auch,g,orthogonal,zu,h]]

\bigskip

Now, the annotation module identifies for each sentence its status: Is it an annotation, an assumption or a claim? Moreover, it extracts the occuring referents and formalizes the content of the respective claim. Thus, a formulation like ``$g\parallel l$ und $l$ ist ausserdem orthogonal zu $h$'' will automatically be formalized as a conjunction ``$g\parallel l \wedge l\perp h$''. Here is the output of the annotation module (with two automatically inserted lines deleted for the sake of brevity):
\bigskip

[[1,[],[],ann,bam,[]],

[2,[g,h,l],[],ang,dkl,[[g,is,line],[h,is,line],[l,is,line]]],

[3,[g,l,h],[],ang,[],[[g,parallel,l],and,[l,orthogonal,h]]],

[4,[g,h],[],beh,[],[g,orthogonal,h]],

[5,[],[],ann,bem,[]]]

\bigskip

As one can see, each sentence is represented as a $6$-tuple 

\begin{center}
(Id,Refs,Names,Status,Function,Content), 
\end{center}

\noindent
where Id is the line number, Refs is a list of referents occuring in the Content, 
Names contains namings of the line for later refrence (``By l1, we have...''), Status tells us whether the sentence is an annotation, an assumption or a claim, Function is a further subdivision of these categories (for example, ``bam'' in line $1$ tells us that this annotation is a start marker for a proof while, ``dkl'' in line $3$ tells us that this assumption is a declaration) and Content is the formalized content (if any) of that line. Note that lines $1$ and $5$ do not correspond to a part of the original text and serve as a starting and an ending marker for the proof (had we started and ended the proof with explicit markers like ``Beweis:'' and ``qed'', it would look the same). The annotation consists of two steps, namely a natural language parser and a formalization routine. For example, for line $2$, the parser output is

\begin{center}
dcl(dclip(es, seien), dcl(sdcl(val(var("var1"), val(var("var2"), val(var("var3")))), type(line(gerade))))).
\end{center}

This tells us that the whole sentence is a declaration (``dcl''), which is started by a declaration initial phrase (``dclip'' and followed by the actual declaration content, which is a simultaneous declaration of several variables (sdcl), each of which receives the type ``line''. From this, the formalization [[g,is,line],[h,is,line],[l,is,line]] is obtained in the second step. 

This list representation is now passed on to the ``text structure'' module; the task of this module is to determine which lines are logically accessible from which other lines. For example, in a proof by case distinction, the case assumption of case $1$ must no longer be available when considering case $2$, so that assumption must not be accessible from the work on case $2$. The same holds true for declarations of variables. This problem of determining the ``range'' of an assumption in a natural language proof is not easy in general; in textbooks and research papers, this often relies on the ability of the reader to infer it from the context from pragmatic considerations. The accessibility rules for Diproche texts are designed to be easy to remember and natural at the same time: In general, an assumption is accessible from all later sentences in the same paragraph, except when it is made in the paragraph that comes directly after a proof starting marker, in which case it is accessible from all later sentences up to the corresponding proof ending marker. Let us look at an example:

\begin{tcolorbox}

(1) Beweis: (2) Es sei $n$ eine nat\"urliche Zahl. 

\medskip

(3) Angenommen, $n$ ist gerade. (4) Dann ist auch $n(n+1)$ gerade. 

\medskip

(5) Angenommen, $n$ ist ungerade. (6) Dann ist $(n+1)$ gerade. (7) Also ist  $n(n+1)$ wiederum gerade.

\medskip

(8) Also ist $n(n+1)$ gerade. 

(9) qed.\footnote{(1) Proof: (2) Let $n$ be a natural number. 

(3) Suppose that $n$ is even. (4) Then $n(n+1)$ is even as well. 

(5) Suppose that $n$ is odd. (6) Then $(n+1)$ is even. (7) Hence, (n+1) is once again even. 

(8) Thus $n(n+1)$ is even.

(9) qed.}

\end{tcolorbox}

Here, the declaration (2), as it comes immediately after the proof starting marker (1), holds for the whole proof, up to the corresponding proof end marker (9). 
The assumption (3) is accessible from the claim (4), but not from any other line in the proof, as the paragraph ends after (4). Likewise, assumption (5) is only accessible from lines $6$ and $7$. The only assumption accessible from the finishing line (8) is thus the declaration (2); however, the implications (implicitly) proved in lines (3)--(4) and (5)--(7) are available, so that the checking is succesful. Here, the accessibility relation would be $\{(2, 4), (2, 3), (2, 4), (2, 5), (2, 6), (2, 7), (2, 8), (3, 4), (5, 6), (5, 7)\}$.
\noindent
In our example case, the accessibility relation generated by the text structure model is simply $\{(2, 3), (2, 4), (3, 4)\}$.

The annotated text, together with the accessibility relation, is now passed on to the module that generates ATP-tasks. For every line that contains a claim (either explictly or implicitly, e.g., by an existential presupposition as in ``Let $n$ be a natural number such that $n^2=n+1$''), the set of accessible assumptions is determined, along with the claims that were already deduced from these assumptions, the implications between earlier claims and the assumptions from which these are supposed to follow and the accessible declarations. The tuple 

\begin{center}
((Assumptions,Claims,FormerTasks,Declarations),Goal)
\end{center}

\noindent
is then passed on to the ATP. In our case, only line $4$ contains a claim; the corresponding ATP-task is this: 

\begin{center}
[[[[g,parallel,l],and,[l,orthogonal,h]]],[],[],[[g,is,line],[h,is,line],[l,is,line]],[g,orthogonal,h]]
\end{center}

The ATP component used for this particular text now automatically breaks the conjunction apart, so that the statements ``$g$ is parallel to $l$'' and ``$l$ is perpendicular to $h$'' are available; also, it has an inference rule that allows one to infer that $g$ is perpendicular to $h$ from the assumptions that $g$ is parallel to $l$ and $l$ is perpendicular to $h$. Thus, the ATP will succeed in verifying the claim based on the available assumptions and the line is checked positively. 

\begin{figure}[h!]
\caption{Entering a Geometric Argument in the Diproche Interface}
\includegraphics[width=15cm]{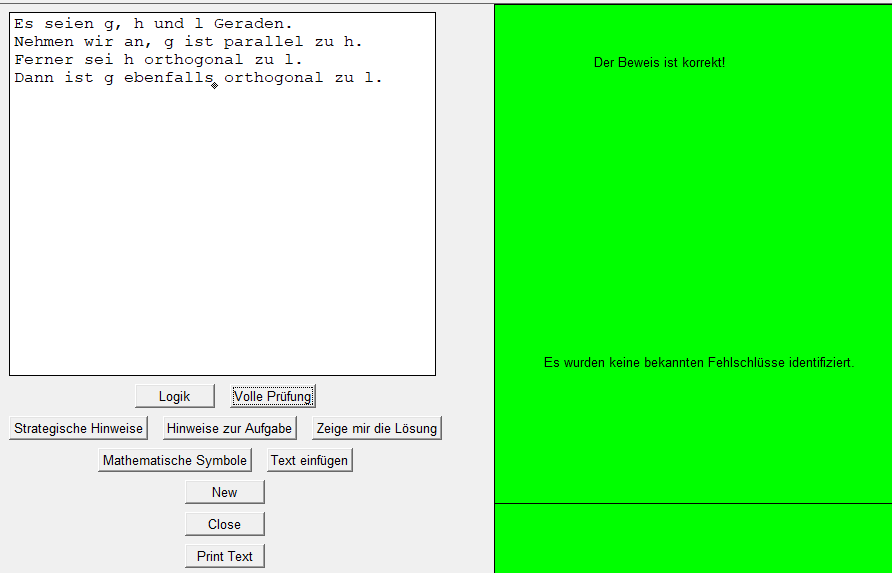}
\end{figure}

\subsection{Example Texts}

In this section, we consider some texts accepted as correct by the current Diproche version.

Our first example is a basic exercise in number theory like those found in Chartrand et al. \cite{CPZ}. It belongs to the ``playing field'' of elementary number theory.


\begin{tcolorbox}
Es sei $n$ eine natuerliche Zahl. Angenommen $n^2$ ist gerade
Wir zeigen: Dann ist $4$ ein Teiler von $n^2$. 

\bigskip

Beweis:

\medskip

Angenommen, $n$ ist ungerade. Dann ist auch $n^2$ ungerade. Widerspruch.

\medskip

Also ist $n$ gerade. Folglich existiert eine natuerliche Zahl $k$ mit $n=2*k$. 
Es sei $k$ eine natuerliche Zahl mit $n=2*k$. Dann folgt $n^2=(2*k)^2=4*k^2$. Also ist $4$ ein Teiler von $n^2$. 
qed.\footnote{Let $n$ be a natural number. Supose that $n^2$ is even. We show: Then $4$ is divisor of $n^2$. 

\smallskip

Proof: 

\smallskip

Suppose that $n$ is odd. Then $n^2$ is also odd. Contradiction. 

\smallskip

Hence $n$ is even. Thus, there is a natural number $k$ such that $n=2*k$. Let $k$ be a natural number with $n=2*k$. Then we have $n^2=(2*k)^2=4*k^2$. Hence $4$ is a divisor of $n^2$. qed.} 

\end{tcolorbox}

A typical example of a solution for a geometry exercise\footnote{The exercise comes from the exercise sheets for the axiomatic geometry lecture at the University of Flensburg by H. Lorenzen.} in Diproche is the following\footnote{Here, $d(a,b,c)$ denotes the triangle with vertices $a$, $b$, $c$, $v(a,b,c,d)$ is the quadrangle with vertices $a$, $b$, $c$, $d$, $s(a,b)$ is the line segment with endpoints $a$ and $b$ and $l(a,b)$ is the line through $a$ and $b$.}:

\begin{tcolorbox}
Es seien $a$, $d$, $c$, $d_1$ Punkte. Es sei $d(a,c,d)$ gleichschenklig. Ferner sei $d(a,c,d)$ rechtwinklig. 
Es sei $l$ der Mittelpunkt von $s(d,d_1)$. 
Angenommen, $l$ liegt auf $l(a,c)$. 
Es sei $l(d,d_1)$ orthogonal zu $l(a,c)$. 

\medskip

Wir zeigen: Dann ist $v(a,d,c,d_1)$ ein Quadrat. 

\medskip

Beweis: Wir haben $s(a,d)\sim s(d,c)$. 
Es ist $l(a,d)$ orthogonal zu $l(c,d)$. 
Es gilt $l$ liegt auf $l(d,d1)$. 
Also gilt $l(l,d)=l(d,d_1)$. 
Damit ist $l(l,d)$ orthogonal zu $l(a,c)$. 
Nach der Mittellotregel folgt $s(a,l)\sim s(l,c)$. 
Also ist $l$ der Mittelpunkt von $s(a,c)$. 
Damit ist $v(a,d,c,d_1)$ ein Parallelogramm. 
Wegen $s(a,d)\sim s(d,c)$ ist $v(a,d,c,d_1)$ sogar eine Raute. 
Also ist $v(a,d,c,d_1)$ ein Quadrat. 
qed.\footnote{Let $a$, $d$, $c$, $d_1$ be points. Suppose that $d(a,b,c)$ is isosceless. Further, let $d(a,b,c)$ be isosceless. Let $l$ be the midpoint of $s(d,d_1)$. 
Suppose that $l$ lies on $l(a,c)$. Let $l(d,d_1)$ be perpendicular to $l(a,c)$. 

\smallskip 

We show: Then $v(a,b,c,d)$ is a square. 

\smallskip 

Proof: We have $s(a,d)\sim s(d,c)$. Let $l(a,d)$ be perpendicular to $l(c,d)$. Then $l$ lies on $l(d,d_1)$. Hence, we have $l(l,d)=l(d,d_1)$. Thus $l(l,d)$ is perpendicular to $l(a,c)$. By the perpendicular bisector rule, it follows that $s(a,l)\sim s(l,c)$. Thus, $l$ is the midpoint of $s(a,c)$. It follows that $v(a,d,c,d_1)$ is a parallelogramm. As $s(a,d)\sim s(d,c)$, $v(a,d,c,d_1)$ is actually a rhombus. Hence $v(a,d,c,d_1)$ is a square. qed.}

\end{tcolorbox}

Finally, we give the following example, which is a version of Thales' theorem in Euclidean geometry\footnote{The proof is the one in Lorenzen \cite{Lo0}, adapted to the Diproche language.}, to indicate the use of annotations for directing the construction of sub-goals in an equivalence proof.

\begin{tcolorbox}
Es sei $d(a,b,c)$ ein echtes Dreieck. 
Es sei $m$ der Mittelpunkt von $s(a,b)$. 

Wir zeigen: Dann ist $d(a,b,c)$ rechtwinklig gdw $s(m,a)\sim s(m,c)$.

\medskip

Beweis: 
Es sei $l:=m(s(a,c))$. 
Dann folgt $l(m,l)||l(b,c)$.  

\medskip

$\Rightarrow$  Es sei d(a,b,c) rechtwinklig. 
Dann ist $l(a,c)$ orthogonal zu $l(b,c)$. 
Also ist $l(m,l)$ orthogonal zu $l(a,c)$. 
Damit folgt $s(m,a)\sim s(m,c)$. 
qed.

$\Leftarrow$ Nun gelte $s(m,a)\sim s(m,c)$. 
Dann ist $l(m,l)$ senkrecht zu $l(a,c)$.  
Also ist $l(b,c)$ orthogonal zu $l(a,c)$.
Damit ist $d(a,b,c)$ rechtwinklig. 
qed.

Also ist $d(a,b,c)$ rechtwinklig gdw $s(m,a)\sim s(m,c)$. 

qed.\footnote{Let $d(a,b,c)$ be a proper triangle. Let $m$ be the midpoint of $s(a,b)$. \smallskip We show: Then $d(a,b,c)$ is right-angled if and only if $s(m,a)\sim s(m,c)$. 

\smallskip

$\Rightarrow$ Suppose that $d(a,b,c)$ is right-angled. Then $l(a,c)$ is perpendicular to $l(b,c)$. Hence, $l(m,l)$ is perpendicular to $l(a,c)$. It follows that $s(m,a)\sim s(m,c)$. qed. 

\smallskip

$\Leftarrow$ Now suppose that $s(m,a)\sim s(m,c)$. Then $l(m,l)$ is perpendicular to $l(a,c)$. Hence $l(b,c)$ is perpendicular to $l(a,c)$. Thus $d(a,b,c)$ is right-angled. qed.

\smallskip

Thus $d(a,b,c)$ is right-angled if and only if $s(m,a)\sim s(m,c)$. 
qed.}
\end{tcolorbox}


\subsection{Further Functions: Goal-Checking, Type-Checking, Hints and Mistake Diagnosis}

The result of the logical check, i.e. the check for the soundness of the inferences of the occuring steps, is not the only kind of feedback that is important to students and it is also not the only kind of feedback that Diproche provides. Without going into detail, we explain here four further kinds of feedback besides logical verifiability that Diproche provides.

\subsubsection{Goal-Checking}

Goal announcements and -modifications are important parts of natural language proofs. For example, one may start by announcing to prove that $A\leftrightarrow B$; then write ``$\Rightarrow$'' to indicate that one is now going to show that $A\rightarrow B$; and then assume $A$, so that the ew goal becomes $B$. Properly mastering such announcements is part 
of learning how to present proofs in natural language. Moreover, in checking an argument, it is not only important whether all steps were sound, but also whether the argument does actually support the consequence it was supposed to prove. In the Diproche system, this is handled by the goal-checker. The goal checker generates a finite (and possibly empty) list of possible goals for each line of the proof and, for each proof end marker, evaluates whether one of the listed goals was achieved and whether it was achieved under the right assumptions. In the example above, if the current goal list consists of the one element $A\rightarrow B$ and $A$ is assumed, then the new goal list is $[B,A\rightarrow B]$; now if the claim $B$ is obtained under no further global assumptions, one can finish this part of the proof by ``qed'' and the goal-checker will accept. If the goal is not reached or if it is reached under additional assumptions, an error message is returned. 
For a detailed explanation and an example of how the goal tracer works, see \cite{CK}.

\subsubsection{Type-Checking}

A mistake that students frequently make is that variables are either not introduced before they are used or that operators are used in the wrong way, e.g., by putting implication arrows between sets, applying set operators to numbers etc. Due to the possibility to enter free text in Diproche, users are free to make such mistakes. Since learning is supported by making mistakes and improving, we regard this as a feature of the system. To be helpful, however, such mistakes in the use of types should be separated from logical mistakes like non-verifiable deduction steps. Therefore, Diproche uses a type-checking algorithm that checks, for each use of a variable, whether the variable has been introduced and whether it is used in accordance with the type that it was assigned when it was introduced. 

\subsubsection{Hints} 

In case a student gets stuck on an exercise, Diproche offers three types of hints: Hints entered manually by the teacher, general hints based on the form of the proof goal or the available assumptions such as ``to prove $A\wedge B$, first prove $A$ and then prove $B$'' and intermediate steps obtained from the completion of the proof by an ATP. The latter two kinds of hints are currently only available for propositional calculus (cf. \cite{CK}). Moreover, the ATP is currently quite weak. We regard it as an interesting, but challenging sub-project to improve this kind of hint; the main obstacle here is of course to generate steps that, in contrast to the output of, say, a Tableau prover, are actually helpful for a beginner writing a ``natural'' proof.



\subsubsection{Mistake Diagnosis} 

In addition to simply marking an inference step as non-verifiable, it is often possible and helpful to identify a particular misunderstanding that caused a fallacious 
inference step. There are various such formal fallacies that occur rather frequently, such as ``deducing'' $\neg b$ from $a\rightarrow b$ and $\neg a$ or ``distributing'' $\neg$ over $\wedge$ to ``deduce'' $\neg a\wedge\neg b$ from $\neg(a\wedge b)$. To provide specific feedback on formal fallacies, Diproche is equipped with a so-called ``anti-ATP'': the anti-ATP works like an ATP, but instead of sound inference rules, it uses common formal fallacies. When the ATP fails to verify an inference step made by the user, the anti-ATP tries to obtain the respective conclusion. If it is successful, the identifier of the rule by which the conclusion was obtained is used to generate a feedback for the user.\footnote{Note that the function of the anti-ATP is neither covered by, nor does it cover, the function of, e.g., a counter-model builder that generates counterexamples when a certain inference cannot be made. While counterexamples are certainly educationally valuable, they do not help in identifying specific and systematic reasoning errors, which is the point of the anti-ATP. On the other hand, the anti-ATP provides no help in situations that do not match with any of the implemented fallacies.} The anti-ATP also has a submodule for identifying typical mistakes in term manipulations, such as distributive use of exponentiation over addition (as in ``$(a+b)^2=a^2+b^2$ '') etc.

For a detailed presentation of the Anti-ATP, see \cite{C1} and the brief discussion in \cite{CK}.

\subsection{Problem Generation}

For several areas, many proving exercises have a common form, which makes it possible to automatically generate exercises. This is desirable as it expands the amount of available training material for the student indefinitely; thus, no matter how many worked-out examples one has seen, one never runs out of ``fresh'' exercises. 

Currently, problem generators have been implemented for the following types of problems:

\begin{itemize}
	\item Propositional Logic: A propositional tautology of bounded length and number of variables is automatically generated; the task is to prove it.
	\item Boolean Set Theory: Two Boolean set terms $t_{0}$, $t_{1}$ (combinations of set variables with unions, intersections and complementation) are generated and the goal is to prove that it holds in general that $t_{0}\subseteq t_{1}$ or that $t_{0}=t_{1}$.
	\item Odd/Even (direct proof) A polynomial $p$ of degree $\leq 3$ with integer coefficients is generated and the goal is to prove that, when $n$ is odd/even, then $p(n)$ is odd/even.
	\item Odd/Even (proof by contraposition) Similarly as for the last type, but the goal is now to show that $n$ is odd/even when $p(n)$ is odd/even.
	\item Odd/Even (proof by case distinction) Similarly as for the last type, but the goal is now to show that $p(n)$ is always odd/even.
	\item Odd/Even (equivalence proof) Two polynomials $p$ and $q$ are generated and the goal is to show that $p(n)$ is odd/even if and only if $q(n)$ is odd/even.
	\item Induction (Divisibility) The goal is to show that, for a fixed natural number $k$, $k$ divides a term of the form $a\cdot b^{cn}+d$, for all natural numbers $n$.
	\item Induction (Inequality) The goal is to show that, for each natural number $n$ larger than a given natural number $k$, a term of the form $a\cdot b^{c\cdot n}+d$ is always less or equal than a term of the form $p\cdot q^{r\cdot n}$.
\end{itemize}

\subsection{Formalization Exercises}

While we will make it clear below that we do not subscribe to a ``formal logic first''-approach to teaching how to prove, a certain mastery of formal language is a necessary prerequisite for writing proofs. 
To support the acquisition of formalization abilities, we implemented two programs for automated formalization exercises: ``Mathe-Diktate'' (``Math Dictations'') and the ``Game of Def''. In ``Math Dictations'', the user is given a mathematical statement like ``the real function $f$ is always larger than the real function $g$'' and is asked to write a statement in the language of first-order logic that formalizes this statement. The users' input formula $\phi$ is checked by using a (strictly ressource-bounded) tableau prover to check whether $\phi\leftrightarrow\psi$, where $\psi$ is the stored template solution. The ``Game of Def'', on the other hand, presents the user with a $21\times 21$-square grid, some squares of which are coloured; the task is to enter a formula in first-order predicate logic that describes the set of coloured squares, using basic predicates like ``is a neighbour of'' or ``is to the right of''. Details on these programs can be found in \cite{C}.

\section{The Language of Diproche}

The language of Diproche is designed to capture a fragment of the German\footnote{A translation to, e.g., English would not be too much effort.} language comprehensive enough to allow natural presentations of solutions to proof exercises in beginner exercises. 

The linguistic units in which such solutions are expressed fall into one of the following categories:

\begin{itemize}
	\item Assumptions (``Suppose that $a$ is parallel to $b$'')
	\item Claims (``Hence $n$ is even''), with justified claims (``Since $n$ is even, there is $k$ such that $n=2k$'') and multiple claims (``Hence $\phi$, so we have, $\psi$, and consequently, we have $\xi$'') as special cases.
	\item Declarations (``Let $n$ be a natural number'') and declarations in connection with a claim (```Let $n$ be a natural number such that $n^2=n$'').
	\item Definitions (``Define $M$ to be the midpoint of $\overline{AB}$''; ``Let $l:=AB$.'')
	\item Goal announcements (``We will show that...''), including subgoal markers like ``$\Rightarrow$'', ``$\Leftarrow$'', ``$\subseteq$'', ``$\supseteq$'' in the proofs of an equivalence or a set equality. 
	\item Start and end markers for proofs (``Proof:'', ``qed'').
	\item Method announcements (``By induction...'', ``By contradiction...'', ``By case distinction...'')
\end{itemize}

For each of these categories, the language contains the usual German triggering expressions. These categories should be self-explaining. The difference 
between assumptions and definitions and declarations (with or without content) is that the latter introduce variables while the former do not. 
Thus, "Assume that $x$ is even" does not serve as a declaration of $x$ and would thus lead to a type mistake if $x$ was not introduced before. 

The following example should illustrate the difference:







(1) "Let $k$ be a natural number such that $n=2k$."

This presupposes that there is such a $k$. Thus, at this point, the checker generates as a subgoal the existence of such a $k$. As it stands, this presupposition would be flagged as unverifable. This would change if it was e.g. preceded by the assumption that $n$ is even.

\medskip

(2) ``Suppose that $n$ is even. Then there is a natural number $k$ such that $n=2k$. Let $k$ be a natural number such that $n=2k$. Then $n^2=(2k)^2=2(2k)$. Thus $n^2$ is even.''

Here, ``$n^2$ is even'' should only be taken to depend on the assumption that $n$ is even, and not on the naming introduced in the third sentence. However, the fourth sentence cleary does depend on this naming as an assumption. In Diproche, this is handled by a selection routine that lists all those namings as assumptions that concern the variables occuring in the respective claim.



\subsection{Formal Expressions}

A typical feature of mathematical language is the mixing of formal expressions with natural language (cf. \cite{CFKKSV}); natural language sentences may contain terms or formulas, or the text may be interrupted by a chain of term manipulations, after which the text continues. Often, such formal expressions come in forms that violate a strict formal syntax, for example in the case of inequality chains like $a=b=1\geq 0$ (strictly speaking, this is not a well-formed expression, as, for example, the second equality sign has propositions, rather than real numbers, on both sides). The formula syntax of Diproche is designed to capture such phenomena by allowing a somewhat ``relaxed'' syntax. In particular, the following expressions can be used in the ``playing field'' of elementary number theory:

\begin{itemize}
\item Arithmetical terms like $a^2+5*(x+2)-3$. When no full bracketing is provided by the user, it is automatically supplemented following the usual priority rules.
\item Inequality chains, i.e., finite sequences alternating between terms and elements of the set $\{=,<,>,\leq,\geq\}$ (where the first and the last elements need to be terms).
\item Manipulation chains, i.e., finite sequences alternating between equalities/inequalities with two sides (no chains) and elements of the set $\{<=>,=>\}$; the (bi)implication sign can also be annotated with a manipulation like $(+3)$ to indicate that the next (in)equality in the chain arises by applying the respective operation. 
Thus, one could, for example, write $a=b\Leftrightarrow(-b)(a-b)=0<=>(*5)(5*a-5*b)=0$.
\end{itemize}

In geometry, we have, for example notations for the segment $s(a,b)$ given by two points $a$ and $b$, the line $l(a,b)$ through two distinct points $a$ and $b$, the triangle $d(a,b,c)$ with vertices $a$, $b$ and $c$,  the quadrangle $v(a,b,c,d)$ with vertices $a$, $b$, $c$, $d$ etc. Moreover, we use $a\in l$ to say that point $a$ lies on line $l$, $g||h$ to say that $g$ is parallel to $h$ etc. All of these have natural language counterparts that one can also use.



\section{Elementary Number Theory and Axiomatic Geometry as Introductory Topics}

In this section, we motivate elementary number theory and axiomatic geometry as introductory topics in learning how to prove and give details about their implementation in Diproche.

\subsection{Criteria for a Suitable Topic for a ``playing field''}


Although proofs occur everywhere in mathematics, not every field is equally suitable as a ``playing field''. In order to be both of use in the teaching of how to prove and work well with the Diproche system, the choice of a playing field depends both on didactical and on technical considerations. We give here some criteria that we used in the determination of suitable playing fields:

\begin{enumerate}
\item A ``flat'' ontology, i.e., a small, fixed number of basic types, rather than a type hierarchy. 
\item A small and fixed language.
\item Proof steps should be reducible to a surveyable set of inference rules (though this set may well contain a few hundred rules, much more than one would want to handle explictly).
\item Ideally axiomatic foundations.\footnote{Though this not a necessary requirement: When the ``common inferential practice'' can still be learned from a corpus study, it can 
	be encoded in inference rules; for example, though number theory is of course axiomatized by the Peano axioms, these play no role in the corresponding ``playing field'' described below.}
\item There should be a rich amount of natural proofs that are ``close'' to their formalization.\footnote{This is, for example, not the case when sophisticated coding machinery is used to formalize finite sequences in the domain of number theory.}
\item The topic should be easy to grasp; its objects should either be familiar to students or one should easily become familiar with them. In this way, the frequent quadruple-difficulty~-- new topic, new level of abstraction, new language, new methodology~-- that beginner students frequently encounter is reduced.
\item There should be clear degrees of difficulty, i.e., clusters of inference rules that allow solutions of many problems.
\item There should be large clusters of ``independent'' problems that do not rely on each other.\footnote{It is quite possible to also implement ``series'' of problems that hierarchically rely on each other, though this has so far not been done. This would correspond to an exercise that is split into several parts. But if the whole field always develops upwards and has no ``levels'' at which one can train, this is a problem.}
\item There should be many exercises that have natural solutions based on self-sufficient text, without reliance on diagrams, intuition etc. 
\end{enumerate}

Elementary number theory satisfies all of these points, while axiomatic geometry satisfies all except possibly (7). In contrast, elementary combinatorics at least fails (1)--(5) and should thus be regarded as a ``bad'' topic for our purposes.\footnote{See, however, Haven \cite{Ha} for an approach to teach stochastics with the help of a automatic system of formal mathematics, namely the interactive theorem prover Coq.}





\subsection{A Very Brief Introduction to Axiomatic Geometry}

In order to keep the paper self-contained, we give here a very brief introduction to the kind of axiomatic geometry that the current Diproche version supports. This is based on the course about axiomatic geometry taught by H. Lorenzen at the EUF in Flensburg as a regular and mandatory part of the curriculum, a course that students usually take in their second semester. The implementation of axiomatic geometry in Diproche is based on the lecture notes by H. Lorenzen \cite{Lo0}, in which all of the material below can be found.

Underlying axiomatic geometry is the notion of an incidence structure (see, e.g., \cite{Lo0}), which is a pair $(\mathcal{P},\mathcal{L})$, where $\mathcal{P}$ is a non-empty set, the elements of which are called ``points'' and a set $\mathcal{L}\subseteq\mathfrak{P}(\mathcal{P})$ of subsets of $\mathcal{P}$, the elements of which are called ``lines''. This is then augmented by relations $\sim$ for ``congruence'', $\parallel$ for parallelism and $\perp$ for orthogonality. One then considers the following axioms (cf. \cite{Lo0}): 

\begin{enumerate}
\item For every two different elements $x,y\in\mathcal{P}$, there is exactly one element of $\mathcal{L}$ that contains both of them (i.e., two points determine exactly one line).
\item There are three elements $P$, $Q$, $R$ in $\mathcal{P}$ such that no element of $\mathcal{L}$ contains all of them (i.e., there are three non-collinear points).
\item For every element $P$ of $\mathcal{P}$ and every element $l$ of $\mathcal{L}$, there is exactly one element in $\mathcal{L}$ which contains $P$ and is either disjoint from $l$ or identical to $l$ (i.e., for every line and every point, there is a unique parallel to the line passing through the point).
\end{enumerate}

An incidence structure satisfying (1)--(3) is called an ``affine plane''; already in this very restricted setting, surprisingly many meaningful exercise problems can be posed. 
However, the possibilities are considerably increased by adding the notions of congruence $\sim$ and orthogonality $\perp$; formally, $\sim$ is a binary relation on line segments, i.e., on pairs of points (thus, a subset of $\mathcal{P}^{2}\times\mathcal{P}^{2}$), while $\perp$ is a binary relation on lines. 

These are characterized by adding the following axioms:

\begin{enumerate}
\setcounter{enumi}{4}
	\item $\sim$ is an equivalence relation; all line segments of the form $\overline{AA}$ are congruent to each other, but none of them is congruent to any line segment $\overline{AB}$ with $A\neq B$.
	\item For each line $l$ and each point $p$, there is exactly one line $h$ such that $l\perp h$ and $p\in h$. 
	\item When $ABCD$ is a parallelogramm, then $\overline{AB}\sim\overline{CD}$ and $\overline{BC}\sim\overline{AD}$.
	\item When $A\neq B$, $C\neq X$ and $\overline{AC}\sim\overline{BC}$, then $CX\perp AB$ if and only if $XA\sim XB$.
	\item There is a rhombus with a midpoint.
\end{enumerate}

Structures in which the axioms (1)--(9) hold are called ``Euclidean plains''; although there are no notions of length or angle measures, let alone areas, it is sufficient to develop a rich fragment of elementary plane geometry, including, e.g., Thales' theorem.

As usual for a mathematical theory, there is, besides the axioms, a rich amount of statements that frequently occur in arguments and can thus be seen as fundamental for the respective area. Among them are the statements of minimal existence (each line contains at least two points; through each point, there are are at least three lines that pass through the point), the compatibility of orthogonality and parallelism (if two lines are both orthogonal to a given line, they are parallel; parallels to lines orthogonal to a given line are also orthogonal to that line), the existence of a fourth point $D$ making $ABCD$ a parallelogramm for every proper triangle $ABC$ etc. An important part of the theory is formed by special types of quadrangles ((symmetric) trapezoids and kites, parallelogramms, rectangles, rhombuses, squares), the fact that these form a lattice under the inclusion relation and their various characterizations, e.g. via their diagonals (for example, ABCD is a parallelogramm if and only if the midpoint of AC coincides with the midpoint of BD). These fundamental statements are present in the Diproche-ATP in the form of inference rules that allow the corresponding deductions, e.g., deducing that ABCD is a parallelogramm from the statement that its diagonals have the same midpoint; usually, one statement is represented by a cluster of several inference rules to account for equivalent, but formally different, formulations. 



\subsection{Didactical Advantages of Axiomatic Geometry and Elementary Number Theory}

Even if the ``real'' task is learning how to prove rather than learning the subject matter the exercises concern, proofs still needs a subject matter~-- teaching proof techniques ``in abstract'', remote from any particular content, for example in the sense of logical calculus is unlikely to be helpful to beginner students, especially those who struggle with proofs. Advanced formal logic is a way to systematize and reflect on proofs after a proof practice has developed, so didactically, it should come after, not before students learn how to prove. We whole-heartedly agree with a Freudenthal quote contained in Wagenschein \cite{Wa}: One cannot organize an area that one does not know. Without external content, proving will look like a symbol game with arbitrary rules. In this section, we will briefly discuss the advantages of teaching how to prove on the basis of elementary number theory and axiomatic geometry. 

First of all, the subject matter should be such that it does not absorb the attention required for aspects of argumentation, such as correctness and critique of arguments, strategies of argumentation and argument discovery and clarity of presentation. Ideally, the subject matter should be familiar to the students. Second, it should contain statements that are simple enough to understand, but not obviously true, and ideally in some way surprising or interesting, so that a desire to prove can arise. Finally, proving should be developed along with techniques of exploration and discovery. This requires a subject matter in which it is possible to get to insights and conjecture by experimenting and observation. All of this is, e.g., hardly the case when working in abstract algebraic structures that were introduced axiomatically.

In elementary number theory, the subject matter, natural numbers, are well-known to students. It also appears to a be a topic that often triggers some curiosity and interest; lectures and seminars on number theory are usually quite popular among students. It contains simple statements that are both surprising and hard to prove to the degree that some conjectures, like the Goldbach conjecture or the prime twin conjecture, can be explained to $5$th-graders still remain undecided, in spite of centuries of effort. Finally, it is accessible to experimental exploration and observation to the degree that ``experimental number theory'' has developed as a branch of mathematics in its own right. Moreover, many interesting problems can be solved by elementary means. Another feature of number theory is the frequent interaction of logical inferences and numerical calculations or algebraic manipulations.



Similar points can be made in the case of geometry: The objects under consideration~-- points, lines, triangles etc.~-- are well-known to the student and the field is full of simple, but surprising statements. Moreover, geometric investigations usually proceed by drawing and observing figures. This interaction between figure and argument is a didactically particularly relevant property of geometry: It is essential to ``see the general'' in the particular figure that one drew, thus learning a way to use intuition in mathematics. Moreover, geometrical proofs teach valuable heuristical lessons: First, one often proceeds by ``intentional changes'', i.e., by ``viewing objects in a new way'', when, for example, a triangle that emerged as a ``by-product'' of a figure suddenly becomes the center of attention; secondly, it is frequently required to  introduce new objects for the sake of an argument, like auxiliary lines. The axiomatic aspect adds to the experienced interaction between picture and text, since it both makes it possible and requires writing texts that, though strongly ``inspired'' by a picture, need to ``stand for themselves'', without reliance on intuition.\footnote{The way in which picture are in fact part of an argument and not mere illustrations is a fascinating topic in its own right. See, e.g., the work by Mumma \cite{Mumma} on the use of diagrams in Euclid's elements.}


Finally, both number theory and geometry are areas from which many central branches of mathematics, such as algebra or analysis, historically developed; thus, acquaintance with these areas forms a basis for a genetical approach to modern mathematics as in Toeplitz \cite{To}. 

\section{Implementing Elementary Number Theory and Axiomatic Geometry in Diproche}

A playing field is only successfully implemented when a large set of exercises can be given for which Diproche will accept, modulo reformulations, all solutions that would be considered as correct and intended at the relevant level of education. Diproche should not limit the possibilities of argumentation, or it should at least do so as little as possible. 

At the start of the implementation of a new playing field is thus the choice of a corpus of exercises and solutions, preferably written by someone with no relation whatsoever to Diproche. One part of this corpus serves as ``training data'', i.e., it helps to identify relevant linguistical and logical phenomena along with typical inferences for the respective area. In both cases, this required adapting the formula parser to the new notation, adapting the text parser both to the new vocabulary and new grammatical phenomena (for example, the implementation of geometry required adding natural language formulations for ternary predicates such as ``$g$ is the parallel to $l$ through $p$'') and writing one or several new submodules of the automated theorem prover for verification. The other part is used as ``test data'': After the implementation is finished based on the ``training examples'', the problems of the intended kind in the corpus should be solvable using Diproche by texts sufficiently similar~-- again, modulo changes in formulation~-- to the given solutions. When new phenomena are observed in this way, a new ``round'' is started. 

In the implementation of the playing fields for number theory and axiomatic geometry in Diproche, this procedure has not been adhered to very strictly; however, as described below, the actual practice in implementing Diproche bears sufficient resemblance to the above-mentioned strategy. 

\subsection{Elementary Number Theory}

In the case of number theory, the implementation was based on a set of exercises used in an introductory Algebra course at the University of Flensburg, along with their template solutions. It was decided to cover the notions of parity, divisibility, residues, squares and cubes, along with equality and inequality. 

A bunch of rules for handling unary and binary predicates (like ``ungerade'' (odd), ``Quadratzahl'' (square number), ``teilt'' (``divides'')) 
was added to the textparser module, which also allows for collective constructions like ``Let a, b, c be even.'', along with a corresponding formalization routine. 
Symbolic expressions for divisibility ($a|b$) and congruences ($a\sim(m)b$) were added to the formula parser.
The ATP for number theory consists of $5$ submodules (not counting the module handling propositional and first-order logic in general):

\begin{itemize}
\item A module with rules for general number theory with $158$ rules.
\item A module with rules for divisibility with $40$ rules.
\item A module with $58$ rules for congruences.
\item A module with special rules for natural numbers (in contrast to integers) with  $52$ rules.
\item Specials rules for term manipulations, equality and inequality. 
\end{itemize}

The resulting system was then ``tested'' with the problems and exercises in chapter 3 of Chartrand, Polimeni and Zhang \cite{CPZ} that concerned proofs involving odd and even numbers.\footnote{To be precise, we considered the examples 3.4, 3.5, 3.6, 3.8, 3.10, 3.11, 3.12, 3.13, 3.14, 3.15, the proofs of which are presented in the book and which were rewritten in the Diproche syntax, along with the exercises 3.8, 3.9, 3.10, 3.16, 3.17, 3.18, 3.19, 3.20, 3.21, 3.26, 3.27, 3.28, 3.29, for which we wrote the solutions ourselves. A few exercises and examples were excluded from the sample in spite of belonging to the topic of odd and even numbers, since they either used fractions or posed problems to explictly given finite sets (like ``For all $x$ in $\{1,2,3\}$...''), both of which is currently not supported in the number theory module.} The result was encouraging: Almost all problems in that section that were of the intended kind could be solved within Diproche with solutions that were quite similar to what a German translation of the English text would look like.\footnote{The reader may want to consider the solutions in \cite{CPZ} to our first example text above, which treated a completely analogous problem.}\footnote{The translation from English to the Diproche CNL was done by the author. Certainly, it would be desirable to see how well external users perform after some introduction to the system. We plan to take up this point in future work.} Particularly interesting are examples 3.19 and 3.20, which contain delibarately flawed proofs to be checked by the reader: In both cases, the proof text could easily be transcribed into the Diproche language and Diproche detected the flaw (though in 3.19, the explicit assertion that $1$ is odd had to be added so that the problematic step to be detected by the reader became the \emph{only} step that Diproche highlighted as non-verifiable).
Only 3.16 and 3.17 were a bit more difficult: 3.17 required replacing phrases like ``$n$ is even'' or ``$n$ is odd'' by ``$2|n$'' and ``$\neg2|n$'', respectively, in certain places. The reason is a difficulty to handle the priority of logical operators in natural language: Thus, Diproche will read ``not $A$ or $B$'' as ``(not $A$) or $B$'', while here, it should be read as ``not ($A$ or $B$)''. To handle this, the respective phrase had to be formalized, resulting in a considerable deviation from the original text. 3.16, which is a nesting of proof strategies (an equivalence proof, the first part of which uses a case distinction) was problematic due to the use of phrases like ``of the same parity'' and ``of different parities'', which are not implemented in the current system\footnote{It would not be difficult to do so and it might be added in a later version}. Here, considerable modifications to the original argument would have been necessary, including adding argumentation steps that somewhat stray away from the actual goal. For this reason, we regard the system's performance for examples 3.16 and 3.17 as failures. Thus, the system was successfully tested in $22$ out of $24$ cases ($23$ out of $25$ if one takes the two different solutions offered in \cite{CPZ} for example 3.14, both of which could be adapted to Diproche, as two different examples), thus yielding a success rate of about $92$ percent. On our office computer, the average running time for these cases was about $7$ seconds, with a maximal running time of about 20 seconds. 

In addition, the first $5$ of the $8$ results in the chapter $4.1$ in \cite{CPZ} on divisibility were successfully reproduced in Diproche. 

We regard this as a quite positive result, especially since it is not the goal of the Diproche system to serve as a general automatic checking device for arbitrary proof exercises, but, much more modest, to provide a tool to practice proving on the basis of didactically suitable exercises. That the system only works for certain types of exercises is fine, as long as this type contains a sufficient amount of didactically suitable material. In this respect, the ability of the system to capture a reasonable amount of exercises from an established textbook is encouraging.\footnote{
To add some anecdotal evidence, our colleague Michael Schmitz from the math department in Flensburg searched the database with German math olympiad problems for proving problems of the desired kind up to 8th grade and was successful in writing solutions accepted by Diproche for two of them, namely MO090833 and MO520833.  Of course, Diproche is not at all designed with math olympiad training in mind.}



\subsection{Axiomatic Geometry}

In the case of geometry, things were somewhat more complicated; as the topic is somewhat special despite its didactical qualities, much less material is available.\footnote{Note that the specific character of the approach to axiomatic geometry pursued here, which does not contain treatments of areas or angles and allows for rather ``wild'' finite, discrete models, the methods often used in automatizing geometrical proving, like the area method, the angle method, algebraic methods etc. are not applicable in this context. Also, databases with problems in elementary plane geometry usually contain little material that is suitable for this approach.} A part of the exercises for the axiomatic geometry course in Flensburg was used to develop the system, and then another part was used for ``testing'', along with statements and proofs in the lecture notes. The axiomatic methodology of the course, due to which it lends itself easily to formalization and automatization, is in another respect a source of a considerable difficulty: Since the course continues to develop notions, terminology and methodology even concerning its most basic concepts, it is hard to come up with ``degrees of difficulty'' (i.e., sets of inference rules), that are well-adapted to a considerable number of exercises. Often, after a test case had been successfully processed, we needed to add a rule trivializing that exercise in order to provide a reasonable framework for the succeeding exercises. It is of course possible to always allow the ``full power'' of the geometry module of the ATP, but this would make the system unsuitable for applications during the course. It is still an open challenge to identify reasonable degrees of difficulty in this area. If this fails, we might be forced to specify a different set of inference rules for each single exercise. While this is certainly possible, it is clearly not the most convenient solution.

A new feature of the geometry ``playing field'' was that the domain is many-sorted (consisting of points, lines, segments...) and that the applicability of inference rules depends on the types of the occuring objects. Moreover, new objects are often not introduced by explicit declarations, but by constructions (``Define $m$ as the point of intersection of $l$ and $g$''). Thus, the geometry-ATP needs to perform type computations. 

The geometry ATP was developed by (i) incorporating inference rules reflecting the use of the axioms or basic lemmata in the lecture notes \cite{Lo0} (ii) incorporating inference rules explictly mentioned in the lecture notes \cite{Lo0} and (iii) testing the resulting system against various simple proofs in the lecture notes and solutions to exercises and supplementing the system when it was necessary. At the time of writing, the geometry ATP contains $544$ topic-specific rules, some of which refer to further subrules which are not counted here.


The first $10$ exercise sheets for the geometry exercise were used as a testing sample, though the separation between cases used for development and test cases was not strictly upheld, due to the small amount of material. The $10$ exercise sheets contained a total of $47$ obligatory exercises\footnote{Some sheets additionally contained extra exercises, which are not counted here.} (many of them containing sub-items, which are not counted as separate exercises). Of these, only $12$ were suitable for a Diproche treatment in that they posed proving problems expressible in the geometrical language currently provided by Diproche. The exercises sorted out as ``unsuitable'' exhibited one of the following traits: 

\begin{itemize}
\item They were meta-problems about models of certain axiomatic theories, such as ``find an affine plane such that...'' or ``show that a finite affine plane of order $n$ contains $n^2$ points''. Such problems are not expressible in a geometrical language (at least not in any natural way), but would require a language talking about structures.
\item They took place in some particular finite model of the axioms, like ``check whether Thales' theorem holds in the following model $M$''. Currently, Diproche does not support working in a specific model; though it would be easy to encode such a model, arguments about these typically take place in some meta-language and make heavy use of ``without loss of generality''-arguments via symmetry etc. Such arguments, though not inaccessible in principle to automated checking (one could, e.g., automatically generate and check the ``symmetric'' cases when a symmetry argument is made), are not supported by the current Diproche version.
\item They contained notions from dynamic geometry, in particular reflections. While such notions can clearly be implemented in Diproche, they are currently not. The reason is that the current interface does not allow a convenient and natural encoding of those. This will be changed in future work, when the interface has been re-designed.
\item They were not proving exercises; rather, the goal was, e.g., to count certain objects, to draw an example of a certain objects, to carry out a certain construction or to organize a shuffled set of given sentences to a sound proof.
\end{itemize}

Many of these $12$ contained ``unsuitable'' sub-exercises (we counted an exercise as ``suitable'' when it contained at least one suitable sub-exercises). In the end, we arrived at $23$ suitable sub-exercises. When considering these exercises, the following difficulties soon became apparent: 

\begin{itemize}
\item Formulations in meta-language: ``Through each point, there are at least three lines in an affine plane''. While it is possible to re-formulate this as 
``Let $p$ be a point. Then there are lines $g$, $h$, $l$ such that $\sim l_1=l_2$, $l_2$, $l_3$ are pairwise distinct and $p\in l_1$ and $p\in l_2$ and $p\in l_3$'', the German equivalent of which can be processed by Diproche\footnote{The precise Diproche formulation is as follows: ``Es sei $p$ ein Punkt. Dann existieren Geraden $l$, $g$, $h$ so, dass $\sim g=h$ und $\sim g=l$ und $\sim h=l$ und $p\in g$ und $p\in h$ und $p\in l$.''}, it is quite cumbersome to pose an exercise in this formulation, let alone write up a solutions in such a 
way that one ends up proving exactly this involved statement so that the goal-checker will regard it as reached. Thus, goal-checking should in some cases be ignored and the exercise should count as solved when three lines are defined and proved to be pairwise different that contain $p$.

\item Trivial cases: Typical examples of degenerate cases in geometry are triangles where two or all three vertices coincide or are collinear. Often, theorems hold for such cases, but the argument needs to be modified (quite often, it trivializes). While it is possible to write this up in the form of an explicit case distinction, it is cumbersome to do so. 
For this reason, assumptions excluding trivial cases were occasionally added.
\end{itemize}

For these reasons, we (i) ignored the goal-declaration and goal-checking and regarded the exercises as solved as soon all ``parts'' of the desired conclusion were obtained (as it would usually be done when correcting an exercise) and (ii) reformulated and simplified exercises by adding extra-assumptions that banned degenerate cases. 

With these modifications, $13$ sub-exercises could be successfully adressed with the current version of the Diproche system, i.e., a bit over $50$ percent. Most of these exercises required some change of the source code in order to go through; many of these were bugfixes (like correcting a misnamed variable in the ATP), a few required adding ATP-rules or variants of ATP-rules already present. In one case, items were added to the vocabulary.  For these exercises, the average running time was sufficiently low (less than $5$ seconds on our office computer). For the remaining exercises, we were unsuccessful for one or several of the following reasons:

\begin{enumerate}
\item While the exercise itself could be formulated within the current conceptual scope of Diproche, the intended solution used concepts or clusters of inference rules that are not supported by the current version of the system.	
\item The automatized checking was ``too precise''. For example, when entering a solution, it became apparent that extra arguments were required showing that, e.g., two lines are distinct when this was ``apparent from the picture''. While this can be regarded as a positive effect of using an automated checking system that considers the text ``in itself'' purely logically, without referering to intuition, this becomes cumbersome as soon as the ``creative'' aspect of an exercise outweights its ``logical'' aspect and long and involved texts are needed to present a correct clever idea.
\end{enumerate}

While (1) is a techical problem that can be overcome by amending the system components, (2) seems to be more of a ``sociological'' problem: As the content of a lecture develops, certain kinds of ``sloppyness'' in proofs become acceptable, as convenience in expressing ideas is traded in for precision and strictness. In textbooks, this is sometimes made explicit by remarks saying that ``such cases will be ignored in the future''; typically, these are ``degenerate cases'', where some set is empty, some number is $0$ etc. Whether Diproche should even try to reflect this part of mathematical practice is unclear to us at the moment: It would, on the one hand, considerably extend the amount of approachable exercises. On the other hand, checking would become less reliable, more resemblant to that of a human tutor who may well overlook the fact that a certain ``uninteresting'' special case has been skipped. In order to decide whether this is a good idea, one should study how a computer's feedback is viewed by the students in comparison with a human feedback, which is a question for human-machine-interaction. In any case, these difficulties are partly due to the fact that real-world proof texts are not purely logical, but also social objects, dependent in their acceptability on the social context (an acceptable step in a proof in the ``Annals of Mathematics'' is often not acceptable as a step in a homework for beginner students) and such contexts are hard, if not impossible to capture within the strict boundaries of a logical system. At this point, the application of automated proof checkers in teaching mathematics thus leads to intricate interactions between logic, sociology and psychology. 









Thus, the picture is much less clear in the case of axiomatic geometry than in the case of elementary number theory above. However, it should be noted that the exercises considered in the evaluation of the elementary number theory playing field all belong to a rather restricted type of exercise, while the geometry problems under consideration formed a substantial variety in content, vocabulary, difficulty and available background. In this light, we draw two consequences from the above results:

\begin{enumerate}
	\item Not surprisingly, Diproche cannot be expected to work well for arbitrary exercise problems in axiomatic geometry as given on actual exercise sheets. While this actually turned out to work in several cases, exercises need to be carefully selected and formulated when Diproche is to be used for geometry exercises.
	\item For carefully selected exercises, Diproche accepts proof text that resemble template solutions reasonably well while adding a layer of formal precision. 
\end{enumerate}

\section{Comparison with Other Systems}

There are quite a few educational softwares supporting the construction of proofs in elementary geometry for beginners; the first we are aware of is G. Holland's GEOLOG/GEOBEWEIS\footnote{The system is still available online under \url{https://web.archive.org/web/20141104052837/http://www.staff.uni-giessen.de/~gcp3/Geolog/geolog.html}; also see \url{http://home.mathematik.uni-freiburg.de/didaktik/material_download/Geometrie_Aufsatz/node10.html}.} which was successfully used in high school teaching in Germany in the $1990$s, see, e.g., Holland \cite{Ho} and Lorenzen \cite{Lo1}. A more recent example is the system QED-Tutrix (see, e.g. Font et al. \cite{FRG}\footnote{Also see the QED-Tutrix Homepage. \url{http://turing.scedu.umontreal.ca/qedx/}.}), which offers an interactive tutor ``Turing'' that gives feedback and hints during the proof development. 
In comparison to these systems, Diproche stands out in allowing a free input in natural language and not requiring explicit mentioning of rules, thus making it possible to work in contexts where the number of rules becomes so vast that users cannot be expected to overview and explictly mentioning all of them (cf. \cite{CK}). This freedom of entering free text is also the freedom to make many mistakes that are impossible in systems like those mentioned, for example the use of undefined expressions. 

A system for proof verification for didactical purposes that also allows free text input by the user is ``Lurch'',\footnote{Available from \url{lurchmath.org}.} see, e.g., Carther and Monks \cite{CM}. However, Lurch does no natural language processing; instead, the user is then required to annotate the text by hand, marking the ``meaningful'' parts of the text (formal expressions) either as ``claims'' or as ``reasons'' (inference rules) or as ``premises'' (assumptions required for the application of the inference rule).
Unmarked portions are ignored by the Lurch checking process. One could say roughly that the user is thus doing the work that the automatic annotation routine does in Diproche. 
Thus, while users may write whatever they want, they then have to get clear about the status of the sentences and text parts they use, which is certainly a good exercise. Still, Lurch does not ``understand'' and correct the natural language formulations; also, steps usually have to be justified by explicitly naming the the inference rule by which they are supposed to work, along with the premises used, which soon becomes unnatural and also infeasible when the number of rules becomes too large (however, Lurch does allow for ``smart rules'' that can be used without mentioning the premises, like ``by propositional logic''). 

Another system working in a controlled language is Elfe (see Broda and Dore \cite{BD}). As one can see from the sample texts in \cite{BD}, Diproche takes more steps towards allowing the user to use natural language. Also, Elfe uses professional ATPs for the verification, while Diproche attempts to model the notion of an acceptable step very precisely and flexibly with dependence on the problem by ATPs specifically written for this purpose.

Finally, we mention Concludio, a system mainly developed by Grewing, \cite{Gr}, which is currently tested at the university of Aachen. In Concludio, natural language proofs can be built up from texts fragments that can be chosen from a menue while terms can be manually entered. It differs from Diproche for example in not allowing free text input and also no problem-specific restrictions of the accepted inference steps.

To the best of our knowledge, the anti-ATP is an original feature of Diproche.\footnote{The anti-ATP in the Concludio system was written after the author had communicated the idea to the Concludio developers as a part of a cooperation of the two projects.}

\section{Conclusion and Further Work}

We hope to have made it plausible that adding natural language processing on top of formal verification tools leads to promising tools in the teaching of elementary proving strategies and proof presentation skills. Our experience so far shows that the Diproche system implemented so far works well for short texts with a simple logical structure. The system has so far been tested on two faculty members in Flensburg and we expect to gain some experience of letting students work with it in the next semester.
For the time being, it is encouraging that example proofs from various areas and sources could be successfully entered and checked by the system with only minor changes. 

Although the Diproche language is a fragment of natural language, some learning is still required in order to write Diproche proofs. If the system is to be used in a lecture, natural mathematical language must to a certain degree become a part of education, at least to the degree that formulations for assumptions, claims etc. are explained. (Our experience so far is that users with experience in writing proofs learned the acceptable fragment very quickly from looking at a few examples.) This might be regarded as a disadvantage; however, it could also considerably lower the difficulties that beginner students often have with expressing proofs, even if they have the right idea.

Clearly, there is no lack of possible extensions; exercises, degrees of difficulty and whole playing fields can be added, grammar and vocabulary can be amended to allow for even more natural language formulations; the analysis of acceptable inference rules for several areas at several levels can be made more systematic and substantiated empirically through a corpus analysis; similarly, the set of formal fallacies recognized by the anti-ATP can and should be systematically developed, etc. 

A particular problem with the Diproche approach arises for longer proof texts: As the ATP modules are designed to reflect the notion of an ``elementary proof step'' for the respective context, they contain a vast number of rules, thus considerably slowing down the verification when the number of available premises increases. One possibility to improve the performance would be to make the verification ``smarter'' by using heuristics~-- possibly obtained through machine learning techniques~-- both for selecting relevant premises (see, e.g., \cite{AHKTU}, \cite{KLTUH}, \cite{CKKS} for work in this direction) and inference rules.

While we regard it neither as realistic (nor as desirable) that a system like Diproche could replace a human corrector in the foreseeable future in applications like advanced exams or even math competitions, we regard our successes so far as a sufficient indicator that such systems can indeed be set up to cover a considerable portion of the proving exercises for basic lectures like beginning linear algebra or analysis. 

\section{Acknowledgements}

We thank three anonymous reviewers for providing detailed feedback that helped to improve the presentation.

\bibliographystyle{eptcs}

\begin{thebibliography}{10}
\providecommand{\bibitemdeclare}[2]{}
\providecommand{\surnamestart}{}
\providecommand{\surnameend}{}
\providecommand{\urlprefix}{Available at }
\providecommand{\url}[1]{\texttt{#1}}
\providecommand{\href}[2]{\texttt{#2}}
\providecommand{\urlalt}[2]{\href{#1}{#2}}
\providecommand{\doi}[1]{doi:\urlalt{http://dx.doi.org/#1}{#1}}
\providecommand{\bibinfo}[2]{#2}

\bibitemdeclare{article}{AHKTU}
\bibitem{AHKTU}
\bibinfo{author}{Jesse \surnamestart Alama\surnameend}, \bibinfo{author}{Tom
  \surnamestart Heskes\surnameend}, \bibinfo{author}{Daniel \surnamestart
  Kühlwein\surnameend}, \bibinfo{author}{Evgeni \surnamestart
  Tsivtsivadze\surnameend} \& \bibinfo{author}{Josef \surnamestart
  Urban\surnameend} (\bibinfo{year}{2011}): \emph{\bibinfo{title}{Premise
  Selection for Mathematics by Corpus Analysis and Kernel Methods}}
\newblock {\sl \bibinfo{journal}{Journal of Automated Reasoning}}
  \bibinfo{volume}{52}, \doi{10.1007/s10817-013-9286-5}.

\bibitemdeclare{inproceedings}{BD}
\bibitem{BD}
\bibinfo{author}{K. \surnamestart Broda\surnameend} \&
  \bibinfo{author}{M. \surnamestart Dore\surnameend}
  (\bibinfo{year}{2019}): \emph{\bibinfo{title}{ Towards Intuitive Reasoning in Axiomatic Geometry}}
\newblock pp. \bibinfo{journal}{P. Quaresma, W. Neupert (eds).: Theorem proving components for Educational software (ThEdu'18) EPTCS 290}, \doi{10.4204/EPTCS.290.4}.


\bibitemdeclare{unpublished}{CK}
\bibitem{CK}
\bibinfo{author}{M.~\surnamestart Carl\surnameend} \&
  \bibinfo{author}{R.~\surnamestart Krapf\surnameend} (\bibinfo{year}{2019}):
  \emph{\bibinfo{title}{Das Diproche-System – ein automatisierter Tutor für
  den Einstieg ins Beweisen}}.

\bibitemdeclare{unpublished}{C1}
\bibitem{C1}
\bibinfo{author}{Merlin \surnamestart Carl\surnameend} (\bibinfo{year}{2020}):
\emph{\bibinfo{title}{Using Automated Theorem Provers for Mistake Diagnosis in the Didactics of Mathematics}}.
\newblock \urlprefix\url{arXiv:2002.05083v1}.

\bibitemdeclare{unpublished}{C}
\bibitem{C}
\bibinfo{author}{Merlin \surnamestart Carl\surnameend} (\bibinfo{year}{2020}):
  \emph{\bibinfo{title}{Automatized Evaluation of Formalization Exercises in
  Mathematics}}.
\newblock \urlprefix\url{http://arXiv:2006.01800v1}.

\bibitemdeclare{inproceedings}{CM}
\bibitem{CM}
\bibinfo{author}{Nathan \surnamestart Carter\surnameend} \&
  \bibinfo{author}{Kenneth \surnamestart Monks\surnameend}
  (\bibinfo{year}{2017}): \emph{\bibinfo{title}{A Web-Based Toolkit for
  Mathematical Word Processing Applications with Semantics}}.
\newblock pp. \bibinfo{journal}{ Proceedings of CICM 2017, eds. Herman Geuvers and Jacques Fleuriot, Springer Lecture Notes in Artificial Intelligence (LNAI 10383)}\bibinfo{pages}{272--291}, \doi{10.1007/978-3-319-62075-6_19}.


\bibitemdeclare{inproceedings}{CFKKSV}
\bibitem{CFKKSV}
\bibinfo{author}{Marcos \surnamestart Cramer\surnameend},
  \bibinfo{author}{Bernhard \surnamestart Fisseni\surnameend},
  \bibinfo{author}{Peter \surnamestart Koepke\surnameend},
  \bibinfo{author}{Daniel \surnamestart Kühlwein\surnameend},
  \bibinfo{author}{Bernhard \surnamestart Schröder\surnameend} \&
  \bibinfo{author}{Jip \surnamestart Veldman\surnameend}
  (\bibinfo{year}{2009}): \emph{\bibinfo{title}{The Naproche Project Controlled
  Natural Language Proof Checking of Mathematical Texts}}.
\bibinfo{journal}{Fuchs N.E. (eds) Controlled Natural Language. CNL 2009. Lecture Notes in Computer Science, vol 5972. Springer, Berlin, Heidelberg}
\newblock pp. \bibinfo{pages}{170--186}, \doi{10.1007/978-3-642-14418-9_11}.


\bibitemdeclare{inproceedings}{CKKS}
\bibitem{CKKS}
\bibinfo{author}{Marcos \surnamestart Cramer\surnameend},
  \bibinfo{author}{Peter \surnamestart Koepke\surnameend},
  \bibinfo{author}{Daniel \surnamestart Kühlwein\surnameend} \&
  \bibinfo{author}{Bernhard \surnamestart Schröder\surnameend}
  (\bibinfo{year}{2010}): \emph{\bibinfo{title}{Premise Selection in the
  Naproche System}}.
\bibinfo{journal}{ Giesl J., H\"ahnle R. (eds) Automated Reasoning. IJCAR 2010. Lecture Notes in Computer Science, vol 6173. Springer, Berlin, Heidelberg}
\newblock \bibinfo{volume}{6173}, pp. \bibinfo{pages}{434--440},
  \doi{10.1007/978-3-642-14203-1_37}.


\bibitemdeclare{book}{De}
\bibitem{De}
\bibinfo{author}{R.~\surnamestart Descartes\surnameend} (\bibinfo{year}{1998}):
  \emph{\bibinfo{title}{D. Cress (translator): R. Descartes. Discourse on
  Method and Meditations on First Philosophy}}.
\newblock \bibinfo{publisher}{Hackett Publishing Company
  Indianapolis/Cambridge}.

\bibitemdeclare{article}{BD}
\bibitem{BD}
\bibinfo{author}{Maximilian \surnamestart Dore\surnameend} \&
  \bibinfo{author}{K.~\surnamestart Broda\surnameend} (\bibinfo{year}{2019}):
  \emph{\bibinfo{title}{Towards Intuitive Reasoning in Axiomatic Geometry}}.
\newblock {\sl \bibinfo{journal}{Electronic Proceedings in Theoretical Computer
  Science}} \bibinfo{volume}{290}, pp. \bibinfo{pages}{38--55},
  \doi{10.4204/EPTCS.290.4}.

\bibitemdeclare{article}{FRG}
\bibitem{FRG}
\bibinfo{author}{Ludovic \surnamestart Font\surnameend},
  \bibinfo{author}{Philippe \surnamestart Richard\surnameend} \&
  \bibinfo{author}{Michel \surnamestart Gagnon\surnameend}
  (\bibinfo{year}{2018}): \emph{\bibinfo{title}{Improving QED-Tutrix by
  Automating the Generation of Proofs}}.
\newblock {\sl \bibinfo{journal}{Electronic Proceedings in Theoretical Computer
  Science}} \bibinfo{volume}{267}, \doi{10.4204/EPTCS.267.3}.

\bibitemdeclare{book}{CPZ}
\bibitem{CPZ}
\bibinfo{author}{P.~Zhang \surnamestart G.~Chartrand\surnameend, A.~Polimeni}
  (\bibinfo{year}{2013}): \emph{\bibinfo{title}{Mathematical Proof - A
  Transition to Advanced Mathematics}}, \bibinfo{edition}{third} edition.
\newblock \bibinfo{publisher}{Pearson}.

\bibitemdeclare{misc}{Gr}
\bibitem{Gr}
\bibinfo{author}{F.~\surnamestart Grewing\surnameend} (\bibinfo{year}{2019}):
  \emph{\bibinfo{title}{Concludio Homepage}}.
\newblock \urlprefix\url{https://www.concludio.education/}.

\bibitemdeclare{mastersthesis}{Ha}
\bibitem{Ha}
\bibinfo{author}{Andrew \surnamestart Haven\surnameend} (\bibinfo{year}{2013}):
  \emph{\bibinfo{title}{Automated proof checking in introductory discrete
  mathematics classes}}
\bibinfo{note}{Master Thesis, MIT}.

\bibitemdeclare{book}{Ho}
\bibitem{Ho}
\bibinfo{author}{G.~\surnamestart Holland\surnameend} (\bibinfo{year}{1996}):
  \emph{\bibinfo{title}{GEOLOG-WIN : Konstruieren, Berechnen, Beweisen,
  Problemlösen mit dem Computer im Geometrie-Unterricht der Sekundarstufe}}.
\newblock \bibinfo{publisher}{D\"ummler Bonn}.

\bibitemdeclare{inproceedings}{KLTUH}
\bibitem{KLTUH}
\bibinfo{author}{Daniel \surnamestart Kühlwein\surnameend},
  \bibinfo{author}{Twan \surnamestart van Laarhoven\surnameend},
  \bibinfo{author}{Evgeni \surnamestart Tsivtsivadze\surnameend},
  \bibinfo{author}{Josef \surnamestart Urban\surnameend} \&
  \bibinfo{author}{Tom \surnamestart Heskes\surnameend} (\bibinfo{year}{2012}):
  \emph{\bibinfo{title}{Overview and Evaluation of Premise Selection Techniques
  for Large Theory Mathematics}}.
\bibinfo{journal}{Gramlich B., Miller D., Sattler U. (eds) Automated Reasoning. IJCAR 2012. Lecture Notes in Computer Science, vol 7364. Springer, Berlin, Heidelberg.}
\newblock pp. \bibinfo{pages}{378--392}, \doi{10.1007/978-3-642-31365-3_30}.

\bibitemdeclare{inproceedings}{Lo1}
\bibitem{Lo1}
\bibinfo{author}{Hinrich \surnamestart Lorenzen\surnameend}
  (\bibinfo{year}{1999}): \emph{\bibinfo{title}{Geolog, Geobeweis und Geokon
  – Erfahrungen und Konzepte zum Unterricht}}.
\newblock In \bibinfo{editor}{T.~Weth \surnamestart W.~Herget\surnameend,
  H.~Weigand}, editor: {\sl \bibinfo{booktitle}{Standardthemen des
  Mathematikunterrichts in moderner Sicht. Bericht \"uber die $17$.
  Arbeitstagung des Arbeitskreises ``Mathematikunterricht und Informatik'' in
  der Gesellschaft f\"ur Didaktik der Mathematik e.V., Wolfenb\"uttel}},
  \bibinfo{publisher}{DIVerlag franzbecker}.

\bibitemdeclare{book}{Lo0}
\bibitem{Lo0}
\bibinfo{author}{Hinrich \surnamestart Lorenzen\surnameend}
  (\bibinfo{year}{2002}): \emph{\bibinfo{title}{Zur Didaktik des begrifflichen
  Denkens in der Geometrieausbildung}}.
\newblock \bibinfo{publisher}{Universit\"at Kiel}.
\bibinfo{note}{Habilitationsschrift}

\bibitemdeclare{article}{Mumma}
\bibitem{Mumma}
\bibinfo{author}{John \surnamestart Mumma\surnameend} (\bibinfo{year}{2010}):
  \emph{\bibinfo{title}{Proofs, pictures, and Euclid}}.
\newblock {\sl \bibinfo{journal}{Synthese}} \bibinfo{volume}{175}, pp.
  \bibinfo{pages}{255--287}, \doi{10.1007/s11229-009-9509-9}.

\bibitemdeclare{book}{To}
\bibitem{To}
\bibinfo{author}{Otto \surnamestart Toeplitz\surnameend}
  (\bibinfo{year}{1949}): \emph{\bibinfo{title}{Die Entwicklung der
  Infinitesimalrechnung. Eine Einleitung in die Infinitesimalrechnung nach der
  genetischen Methode. Erster Band. Aus dem Nachlass herausgegeben von G.
  Köthe}}.
\bibinfo{publisher}{Springer-Verlag Berlin Heidelberg}, \doi{10.1007/978-3-642-49782-7}

\bibitemdeclare{article}{Wa}
\bibitem{Wa}
\bibinfo{author}{Martin \surnamestart Wagenschein\surnameend}
  (\bibinfo{year}{1966}): \emph{\bibinfo{title}{Zum Problem des Genetischen
  Lehrens}}.
\newblock {\sl \bibinfo{journal}{Zeitschrift für Pädagogik}}
  \bibinfo{volume}{12}, pp. \bibinfo{pages}{305--330}.

\end{thebibliography}

\end{document}